\documentclass[a4paper, 12pt]{amsart}
\usepackage{amsfonts, amsmath, amsthm}
\usepackage{amssymb}
\usepackage{latexsym}
\usepackage[dvips]{graphicx}
\usepackage{color}

\newtheorem{thm}{Theorem}[section]

\newtheorem{lem}{Lemma}[section]

\newtheorem{rmk}{Remark}[section]

\numberwithin{equation}{section}
\newcommand{\Z}{\mathbb{Z}}

\newcommand{\R}{\mathbb{R}}
\newcommand{\C}{\mathbb{C}}

\newcommand{\pa}{\partial}
\newcommand{\eps}{\varepsilon}
\newcommand{\jb}[1]{\langle #1 \rangle}
\newcommand{\Jb}[1]{\bigl\langle #1 \bigr\rangle}

\newcommand{\op}[1]{\mathcal{#1}}
\newcommand{\cc}[1]{\overline{#1}}

\DeclareMathOperator{\realpart}{\rm Re}
\DeclareMathOperator{\imagpart}{\rm Im}

\newcommand{\sh}[1]{{#1}^{\sharp}}
\newcommand{\tm}{\tilde{m}}

\newcommand{\bfu}{\mathbf{u}}

\title[DNLS systems with multiple masses]{
Remarks on derivative nonlinear Schr\"odinger systems with multiple masses
}
\author[C.~Li]
{Chunhua Li}

\address{%{\rm Chunhua Li}\\
         Department of Mathematics\\ 
         College of Science\\ 
         Yanbian University \\
         No.977 Gongyuan Road\\ 
         Yanji City, Jilin Province 133002\\
         China}

 \email{sxlch@ybu.edu.cn}
\author[H.~Sunagawa]
{Hideaki Sunagawa}
\address{%{\rm Hideaki Sunagawa}\\
         Department of Mathematics\\ 
         Graduate School of Science\\ 
         Osaka University \\
         Toyonaka, Osaka 560-0043\\
         Japan}

 \email{sunagawa@math.sci.osaka-u.ac.jp}
\date{\today }

\begin{document}

%---------------
\begin{abstract}
We prove global existence of small solutions to the initial value 
problem for a class of cubic derivative nonlinear Schr\"odinger systems 
with the masses satisfying suitable non-resonance relations. 
The large-time asymptotics of the solutions are also shown.
This work is intended to provide a counterpart of the previous paper 
\cite{LS} in which the mass resonance case was treated.
\end{abstract}
%---------------

\maketitle

%{\footnotesize \em 
%Dedicated to Professor Nakao Hayashi on the occasion of
%his sixtieth birthday}

%-------------------------------------------------------------------------%
\section{Introduction} \label{Sec_intro}
%-------------------------------------------------------------------------%
This paper is a sequel to \cite{LS}. We continue the  study of 
the initial value problem
\begin{align}
\left\{\begin{array}{cl}
\mathcal{L}_{m_j} u_j=F_j(\bfu,\pa_x \bfu), & t>0,\ x\in \R,\ j=1,\ldots, N,\\
u_j(0,x)=\varphi_j(x), & x \in \R,\ j=1,\ldots, N,
\end{array}\right.
\label{nls_N}
\end{align}
where $\mathcal{L}_{m_{j}}=i\pa_t +\frac{1}{2m_{j}}\pa_x^2$, 
$i=\sqrt{-1}$, $m_j \in \R\backslash \{0\}$, and 
$\bfu=(u_j(t,x))_{1\le j \le N}$ is a $\C^N$-valued unknown function. 
The nonlinear term $\mathbf{F}=(F_j)_{1\le j\le N}$ is assumed to be 
a cubic homogeneous polynomial in 
$(\bfu,\pa_x \bfu, \overline{\bfu}, \overline{\pa_x \bfu})$. 
Our main interest is how the combinations of $(m_j)_{1\le j \le N}$ and 
the structures of $(F_j)_{1\le j \le N}$ affect the behavior of 
$\bfu(t)$ as $t \to \infty$. 

Let us first recall backgrounds briefly. As is well known, cubic 
nonlinearity is critical when we consider the large-time behavior 
of solutions to nonlinear Schr\"odinger equation. 
In general, cubic nonlinearity must be regarded as 
a long-range perturbation. According to Hayashi-Naumkin \cite{HN}, 
the solution $u(t,x)$ to 
\begin{align}
 i\pa_t u + \frac{1}{2}\pa_x^2 u = |u|^2 u, \qquad
 t>0,\ x\in \R
 \label{nls_1}
\end{align}
behaves like 
\[
 u(t,x)=\frac{1}{\sqrt{it}} \alpha(x/t) 
 e^{i\{\frac{x^2}{2t}  -  |\alpha(x/t)|^2 \log t \}}
 +o(t^{-1/2})
\quad \mbox{in} \ \ L^{\infty}(\R_x)
\]
as $t\to \infty$, where  $\alpha$ is a suitable $\C$-valued function. 
An important consequence of this asymptotic expression is that the solution to 
\eqref{nls_1} decays like $O(t^{-1/2})$ uniformly in $x\in \R$, 
while it does not behave like the free solution. 
In other words, the additional logarithmic factor in the phase 
reflects a typical long-range character of the cubic nonlinear Schr\"odinger 
equations in one space dimension. There are several extensions of this result. 
If we restrict our attentions to the single case ($N=1$), the most general 
case is considered in \cite{HN4} 
(see also the introduction of \cite{SS} for a survey of recent development 
on cubic derivative nonlinear Schr\"odinger equations in one space 
dimension). 
Next let us turn our attentions to the system case ($N\ge 2$). 
Recently, a lot of efforts have been made for the study on 
systems of nonlinear Schr\"odinger equations with multiple masses 
(see e.g., \cite{CC}, \cite{HLN1}, \cite{HLO}, \cite{L}, \cite{HOT}, 
\cite{OS}, \cite{KLS}, \cite{IKS}, \cite{Hir}, etc.). 
An interesting feature in the system case is that the behavior of solutions 
are affected by the combinations of the masses as well as the structure of 
the nonlinearity. Note that similar phenomena can be observed in critical 
nonlinear Klein-Gordon systems 
(see e.g., \cite{Tsu}, \cite{Su1}, \cite{DFX}, \cite{Su2}, \cite{KawSu}, 
\cite{KOS}).  
In the previous paper \cite{LS}, we have considered the case where the masses 
satisfy suitable resonance relations. Roughly speaking, that is the case 
where the operator $\op{J}_m=x+\frac{it}{m}\pa_x$ works well 
through the Leibniz-type  rule 
\begin{align}
 \op{J}_m(f_1f_2f_3) 
 = \frac{\mu_1}{m}(\op{J}_{\mu_1}f_1)f_2f_3
  + \frac{\mu_2}{m}f_1(\op{J}_{\mu_2}f_2)f_3
  + \frac{\mu_3}{m}f_1f_2(\op{J}_{\mu_3}f_3)
\label{J_Leibniz}
\end{align}
which is valid only when $m=\mu_1+\mu_2 +\mu_3$ 
(see the condition (a) in \cite{LS} for a precise expression of the mass 
resonance condition). 
In \cite{LS}, several structural conditions on the cubic nonlinearity 
have been introduced under which the small data global existence holds, 
and time-decay properties of the global solutions have been investigated. When we restrict our attentions to the  2-component system
\begin{align}
\left\{\begin{array}{l}
 \op{L}_{m} u = \kappa \cc{u}^2\pa_x v, \\
 \op{L}_{\mu} v =\lambda u^2 \pa_x u, 
 \end{array}\right.
\qquad 
  t>0,\ x\in \R 
\label{nls_2}
\end{align}
with $m$, $\mu \in \R\backslash\{0\}$ and 
$\kappa, \lambda \in \C\backslash\{0\}$, 
we can see that the conditions of \cite{LS} are 
satisfied only if $\mu=3m$ and $\realpart(\kappa \lambda)<0$, 
$\imagpart(\kappa \lambda)=0$ 
(see Remark~\ref{rem_resonance} below for the detail). 
We also note that the system 
\begin{align}
\left\{
 \begin{array}{l}
 \op{L}_{m} u = \kappa \pa_x (v^3) \\
 \op{L}_{\mu} v =\lambda \pa_x (u^3) 
 \end{array}\right.
 \qquad 
  t>0,\ x\in \R 
\label{nls_3}
\end{align}
is not covered by the result of \cite{LS} regardless of the ratio of $m$ 
and $\mu$.

In the present paper, we are interested in the non-resonance case, i.e., 
the case where \eqref{J_Leibniz} is not valid. 
For simplicity of exposition, we concentrate our attentions to 
the simple model systems \eqref{nls_2} and \eqref{nls_3} with the 
initial condition 
\begin{align}
 u(0,x)=u_0(x), \ \  v(0,x)=v_0(x), \quad x \in \R.
\label{data}
\end{align}
The general $N$-component system \eqref{nls_N} will be discussed 
in the final section.

To state our results, let us introduce some function spaces. 
For $s, \sigma \in \Z_{\ge 0}$, 
we denote by $H^s$ the $L^2$-based Sobolev space of order $s$, 
and the weighted Sobolev space $H^{s,\sigma}$ is defined by 
$\{\phi \in L^2\, |\, \jb{x}^{\sigma} \phi \in H^s \}$, equipped with 
the norm $\|\phi\|_{H^{s,\sigma}}=\|\jb{x}^{\sigma} \phi\|_{H^s}$. 
The main results are as follows: 

%-------------------------
\begin{thm} \label{thm_main}
Assume $m\ne 3\mu$ and $\mu\ne 3m$. 
Let 
$u_0$, $v_0 \in H^2\cap H^{1,1}$, 
and assume $\eps:=\|u_0\|_{H^2\cap H^{1,1}}+\|v_0\|_{H^2\cap H^{1,1}}$ 
is sufficiently 
small. Then the initial value problem \eqref{nls_3}--\eqref{data} 
admits a unique pair of global solutions 
$u$, $v \in C([0,\infty); H^2\cap H^{1,1})$. 
Moreover, there exist $u_+$, $v_+ \in L^2(\R_{x})$ with 
$\hat{u}_+$, $\hat{v}_+ \in L^{\infty}(\R_{\xi})$ such that 
\[
 u(t)
 =
 e^{i\frac{t}{2m}\pa_x^2} u_+ + O(t^{-1/4+\delta})
 \quad \mbox{in}\ L^2(\R_x)
\]
\[
 v(t)
 =
 e^{i\frac{t}{2\mu}\pa_x^2} v_+ + O(t^{-1/4+\delta})
 \quad \mbox{in}\ L^2(\R_x)
\]
and 
\[
 u(t,x)=\sqrt{\frac{m}{i t}}\, \hat{u}_+
 \left( \frac{m x}{t}\right) e^{i\frac{m x^2}{2t}} 
 + O(t^{-3/4+ \delta}),
 \quad \mbox{in}\ L^{\infty}(\R_x)
\]
\[
 v(t,x)=\sqrt{\frac{\mu}{i t}}\, \hat{v}_+
 \left( \frac{\mu x}{t}\right) e^{i\frac{\mu x^2}{2t}} 
 + O(t^{-3/4+ \delta})
 \quad \mbox{in}\ L^{\infty}(\R_x)
\]
as $t \to \infty$, 
where $\delta>0$ can be taken arbitrarily small,
and $\hat{\phi}$ denotes the Fourier transform of $\phi$, i.e., 
\[
 \hat{\phi}(\xi)=\frac{1}{\sqrt{2\pi}} \int_{\R} e^{-iy\xi} \phi(y) dy.
\]
\end{thm}
%-------------------------

%-------------------------
\begin{thm} \label{thm_second}
Assume $\mu \ne 3m$ and $\mu\ne 2m$. 
Then the same assertion as Theorem~\ref{thm_main} holds for 
the initial value problem \eqref{nls_2}--\eqref{data}.
\end{thm}
%-------------------------

%---------------------------
\begin{rmk} 
In Theorem~\ref{thm_second}, we do not need any 
restrictions on $\kappa$ and $\lambda$. 
This should be contrasted with the resonance case $\mu=3m$. 
Indeed, if $\mu=3m$ and $(\kappa,\lambda)=(0,1)$, 
we can show that 
\[
 \lim_{t\to \infty} \|v(t)\|_{L^2}=\infty
\]
with a suitable choice of $(u_0,v_0)$, which implies that $v$ cannot be 
asymptotically free. 
\end{rmk}
%-----------------------------

We close the introduction with the contents of this paper. 
The next section is devoted to preliminaries. 
In Section~\ref{Sec_apriori}, we will get an a priori estimate for 
the solution to \eqref{nls_3}--\eqref{data}, and 
Theorem~\ref{thm_main} will be proved in Section~\ref{Sec_proof1}. 
In Section~\ref{Sec_proof2}, we will give an outline of the proof of 
Theorem~\ref{thm_second}. 
In Section~\ref{Sec_gene}, we will discuss a generalization of 
Theorems~\ref{thm_main} and \ref{thm_second} to the $N$-component system 
\eqref{nls_N}.

%-------------------------------------------------------------------------%
\section{Preliminaries}\label{Sec_prelim}
%-------------------------------------------------------------------------%
In this section, we summarize basic facts related to the Schr\"odinger 
operator $\op{L}_m=i\pa_t+\frac{1}{2m}\pa_x^2$. 
In what follows, several positive constants are denoted by the same letter 
$C$, which may vary from one line to another.

%------------------------
\subsection{The operators $\op{J}_m$ and $\op{P}$}
%------------------------
We set $\op{J}_m=x+\frac{it}{m}\pa_x$ and $\op{P}=x\pa_x +2t\pa_t$.  
As is well-known, these operators have good compatibility with $\op{L}_m$. 
We can check immediately that 
\begin{align}
 [\op{L}_m, \op{J}_m]=0,\quad
 [\op{L}_m,\op{P}]=2\op{L}_m, \quad
 [\pa_x, \op{J}_m]=1,\quad
 [\pa_x, \op{P}]=\pa_x,
\label{com_rel}
\end{align}
where $[\cdot, \cdot]$ stands for the commutator of two linear operators. 
Another important relation is 
\begin{align}
 \op{J}_m\pa_x =\op{P} +2it\op{L}_m,
 \label{id_J_and_P}
\end{align}
which will be used effectively in Section \ref{Sec_apriori}.

%------------------------
\subsection{Factorization of the free evolution group $\op{U}_m(t)$}
%------------------------
We set 
\[
 \bigl( \mathcal{U}_m(t) \phi \bigr)(x)
 :=
 e^{i\frac{t}{2m}\pa_x^2}\phi(x)
 =
 \sqrt{\frac{|m|}{2\pi t}}e^{-i \frac{\pi}{4} \mathrm{sgn}(m)} 
 \int_{\R} e^{im\frac{(x-y)^2}{2t}} \phi(y)dy
\]
for $m \in \R\backslash \{0\}$ and $t>0$. 
We also introduce the scaled Fourier transform $\mathcal{F}_m$ by 
\[
 \bigl( \mathcal{F}_m \phi \bigr)(\xi)
 :=
 |m|^{1/2} e^{-i \frac{\pi}{4} \mathrm{sgn}(m)}\, \hat{\phi}(m\xi)
 =
 \sqrt{\frac{|m|}{2\pi}} e^{-i \frac{\pi}{4} \mathrm{sgn}(m)} 
 \int_{\R} e^{-imy\xi} \phi(y)dy,
\]
as well as auxiliary operators 
\begin{align*}
 &\bigl( \mathcal{M}_m(t) \phi \bigr)(x)
 :=
 e^{im\frac{x^2}{2t}}\phi(x),
 \\
 &\bigl( \mathcal{D}(t)\phi \bigr)(x) 
 := 
 \frac{1}{\sqrt{t}} \phi \left(\frac{x}{t} \right),
 \\ 
  &\mathcal{W}_m(t)  \phi 
 :=
 \mathcal{F}_m \mathcal{M}_m (t) \mathcal{F}_m^{-1} \phi,
\end{align*}
so that $\mathcal{U}_m$ can be decomposed into 
$\mathcal{U}_m=\mathcal{M}_m \mathcal{D} \mathcal{F}_m \mathcal{M}_m
=\mathcal{M}_m \mathcal{D} \mathcal{W}_m \mathcal{F}_m$. 
The following lemmas are well-known. 

%------------

%-------------------
\begin{lem} \label{lemma_asympt}
Let $m$ be a non-zero real constant. We have 
\begin{align*}
 \|
  \phi - \mathcal{M}_m \mathcal{D} \mathcal{F}_m \mathcal{U}_m^{-1}\phi
 \|_{L^{\infty}}
 \le
 Ct^{-3/4} \bigl(\|\phi\|_{L^2} + \|\mathcal{J}_m \phi\|_{L^2} \bigr)
\end{align*}
and
\[
 \|\phi\|_{L^{\infty}} 
 \le 
 t^{-1/2}\|\mathcal{F}_m \mathcal{U}_m^{-1}\phi\|_{L^{\infty}}
 +Ct^{-3/4} (\|\phi\|_{L^2} + \|\mathcal{J}_m \phi\|_{L^2})
\]
for $t\ge 1$. 
\end{lem}
%------------------

%-------------------
\begin{lem} \label{lemma_prod}
Let $m$ be a non-zero real constant. We have
\[
\|\mathcal{F}_m \mathcal{U}_m^{-1} (f_1 f_2 f_3)\|_{L^{\infty}}
\le 
 C 
 \|f_1\|_{L^2} \|f_2 \|_{L^2} \|f_3\|_{L^{\infty}}.
\]
\end{lem}
%------------------

We skip the proof of these lemmas 
(see e.g., \S 3 of \cite{LS} and its references).

%------------------------
\subsection{Smoothing properties}
%------------------------
We recall smoothing properties of the linear Schr\"odinger 
equations. Let $\op{H}$ be the Hilbert transform, that is, 
\[
 \mathcal{H}\psi(x):=
 \frac{1}{\pi} \, \mathrm{p.v.}\int_{\R} \frac{\psi(y)}{x-y}dy.
\]
With a non-negative weight function $\Phi(x)$ and a non-zero real constant 
$m$, let us also define the operator $S_{\Phi,m}$ by
\begin{align*}
 S_{\Phi,m}\psi (x)
:=&
 \left\{ \cosh\biggl(\int_{-\infty}^{x} \Phi(y) dy\biggr) \right\}
 \psi(x)\\
 &- i\, \mathrm{sgn}(m)
 \left\{\sinh\biggl(\int_{-\infty}^{x} \Phi(y) dy\biggr) \right\}
 \mathcal{H}\psi(x).
\end{align*}
Note that $S_{\Phi,m}$ is $L^2$-automorphism and that both 
$\|S_{\Phi,m}\|_{L^2\to L^2}$, $\|S_{\Phi,m}^{-1}\|_{L^2\to L^2}$ 
are dominated by $C\exp(\|\Phi\|_{L^1})$. 

The following two lemmas are used effectively in \S \ref{subsec_L2} 
to overcome the derivative loss coming from the nonlinear term:

%-----------------------
\begin{lem} \label{lemma_smoothing}
Let $m$, $\mu_1$, $\mu_2$ be non-zero real constants. 
Let $f$ be a $\C$-valued smooth function of $(t,x)$, and 
let $\mathbf{w}=(w_1,w_2)$ be a $\C^2$-valued smooth function of $(t,x)$. 
We set $\Phi=\eta |\mathbf{w}|^{2}$ with $\eta\ge 1$, and 
$S=S_{\Phi(t,\cdot),m}$. 
Then we have 
\begin{align*}
 \frac{d}{dt}\|Sf(t)\|_{L^2}^2 
 +&
 \frac{1}{|m|}\int_{\R} \Phi(t,x) \Bigl| S|\pa_x|^{1/2}  f(t,x) \Bigr|^2 dx\\
 &\le
 2 \Bigl| \Jb{Sf(t), S\mathcal{L}_m f(t)}_{L^2}\Bigr| 
 + CB(t) \|f(t)\|_{L^2}^2,
\end{align*}
where
\begin{align*}
 B(t)
 =
 e^{C\eta \|\mathbf{w}\|_{L^2}^2}
 \left\{
  \eta \|\mathbf{w}(t)\|_{W^{1,\infty}}^2 
  + 
  \eta^{3} \|\mathbf{w}(t)\|_{L^{\infty}}^6
  +
  \eta  \sum_{k =1}^{2}\|w_k(t)\|_{L^2} 
  \|\mathcal{L}_{\mu_k} w_k(t)\|_{L^2}
 \right\}
\end{align*}
and the constant $C$ is independent of $\eta$. 
We denote by $W^{s,\infty}$ the $L^{\infty}$-based Sobolev space of 
order $s \in \Z_{\ge 0}$.
\end{lem} 
%-----------------------

%-----------------------
\begin{lem} \label{lemma_aux}
Let $m_1$, $m_2$ be non-zero real constants. 
Let $\mathbf{f}=(f_1,f_2)$,  $\mathbf{w}=(w_1,w_2)$ be $\C^2$-valued 
smooth functions of $x \in \R$. 
Suppose that $q$ is a quadratic homogeneous polynomial in $\mathbf{w}$. 
We set $\Phi=\eta |\mathbf{w}|^{2}$ with $\eta\ge 1$, and 
$S_j=S_{\Phi, m_j}$ for $j =1,2$. Then we have
\begin{align*}
 \left|\Jb{S_1 f_1, S_1 \bigl( q\pa_x f_2 \bigr)}_{L^2} \right|
 &\le
 \frac{C}{\eta} e^{C\eta \|\mathbf{w}\|_{L^2}^2} 
  \sum_{k=1}^{2} 
  \int_{\R} \Phi(x) \Bigl| S_k |\pa_x|^{1/2} f_k(x) \Bigr|^2 dx\\
 &\hspace{4mm}+
 Ce^{C \eta \|\mathbf{w}\|_{L^2}^2} 
 \bigl( 1+ \eta^{2}\|\mathbf{w}\|_{L^2}^4
        +\eta^{2}\|\mathbf{w}\|_{L^{\infty}}^4 \bigr) 
 \|\mathbf{w}\|_{W^{1,\infty}}^2 \|\mathbf{f}\|_{L^2}^2,
\end{align*}
where the constant $C$ is independent of $\eta$.
\end{lem} 
%-----------------------

For the proof, see \S 2 in \cite{HNP} as well as the appendix of \cite{LS}.

%---------------------------------------------------------------------------%
\section{A priori estimate for \eqref{nls_3}--\eqref{data}} \label{Sec_apriori}
%---------------------------------------------------------------------------%
This section is devoted to getting an a priori estimate for the solution to 
\eqref{nls_3}--\eqref{data}. Throughout this section, we always assume that 
$\mu\ne 3m$, $m\ne 3\mu$ and $u_0$, $v_0 \in H^2\cap H^{1,1}$ with 
$\eps=\|u_0\|_{H^2\cap H^{1,1}}+\|v_0\|_{H^2\cap H^{1,1}}$. 
Let $T\in (0,\infty]$, and let $u$, $v \in C([0,T);H^2\cap H^{1,1})$ be 
a pair of solutions to \eqref{nls_3}--\eqref{data} for $t\in [0,T)$. 
We set 
\begin{align*}
 \alpha(t,\xi)=
 \mathcal{F}_{m} \Bigl[ \mathcal{U}_{m}^{-1} u(t,\cdot) \Bigr](\xi), \quad
 \beta(t,\xi)=
 \mathcal{F}_{\mu} \Bigl[ \mathcal{U}_{\mu}^{-1} v(t,\cdot) \Bigr](\xi). 
\end{align*}
We also define
\begin{align*}
 E(T) =& \sup_{(t,\xi) \in [0,T)\times \R} 
  \Bigl[\jb{\xi} \bigl(|\alpha(t,\xi)| + |\beta(t,\xi)| \bigr) \Bigr]\\
&+
\sup_{0\le t< T}
 \biggl[
 (1+t)^{-\frac{\gamma}{3}} \Bigl(
  \|u(t)\|_{H^{2}}+  \|v(t)\|_{H^{2}}  
 + \|\op{J}_{m} u(t)\|_{H^{1}} + \|\op{J}_{\mu} v(t)\|_{H^{1}}
 \Bigr) \biggr]
\end{align*}
with $\gamma \in (0,1/4)$. 
We are going to show the following: 

%---------------
\begin{lem} \label{lemma_apriori}
Assume $m\ne 3\mu$ and $\mu \ne 3m$. 
There exist positive constants 
$\eps_{1}$ and $K$, not depending on T, such that 
\begin{align}
 E(T)\le \eps^{2/3}
 \label{est_before}
\end{align}
implies  
\begin{align*}
 E(T)\le K\eps,
 %\label{est_after}
\end{align*}
provided that $\eps \le \eps_{1}$.
\end{lem}
%----------------

To prove this lemma, it is sufficient to show
\begin{align}
\|u(t)\|_{H^{2}}+  \|v(t)\|_{H^{2}}  
 + \|\op{J}_{m} u(t)\|_{H^{1}} + \|\op{J}_{\mu} v(t)\|_{H^{1}}
\le K_1 \eps (1+t)^{\frac{\gamma}{3}}
\label{est_L2}
\end{align}
and
\begin{align}
|\alpha(t,\xi)| + |\beta(t,\xi)| \le K_2 \eps \jb{\xi}^{-1}
\label{est_ptwise}
\end{align}
under the assumption \eqref{est_before}, where $K_1$ and $K_2$ are positive 
constants not depending on $T$. Then the desired estimate 
follows by setting $K=K_1+K_2$. So our task is to show these two estimates. 
Before going into the proof of them, let us point out the differences between 
the approach in \cite{LS} and the present one. 
In the present setting, we cannot use the operator $\op{J}_m$ directly 
because the Leibniz-type rule is not valid for \eqref{nls_3} without 
growth in $t$. 
Instead of doing so, we use the dilation operator $\op{P}$ through the 
relation \eqref{id_J_and_P}. The facts that $\op{P}$ is independent of 
$m$, $\mu$ and obeys the usual Leibniz rule allow us to obtain the 
first estimate \eqref{est_L2}. 
This idea is originated by Hayashi-Naumkin \cite{HN1} 
in the study of the generalized Benjamin-Ono equation, and applied to 
single cubic derivative nonlinear Schr\"odinger equations in \cite{HN2}, 
\cite{HN3} (see also \cite{SS}). In the derivation of the second estimate 
\eqref{est_ptwise}, we use the factorization of $\op{U}_m(t)$ and reduce the 
original system to simpler equations satisfied by $\alpha$, $\beta$. 
The oscillating factor in the reduced equation enables us to get the 
desired estimate \eqref{est_ptwise} under the non-resonance condition 
$\mu\ne 3m$ and $m\ne 3\mu$. 
Similar idea was used previously in \cite{DFX}, \cite{Su2}, \cite{KawSu}, etc., for the Klein-Gordon case and in \cite{HLN2} for the final state problem for 
a quadratic nonlinear Schr\"odinger system.

%-------------------------------%
\subsection{$L^2$-estimates} \label{subsec_L2}
%--------------------------------%
The goal of this part is to get \eqref{est_L2}. 
Proof is divided into five parts: we will derive four kinds of 
$L^2$-estimates first, and then piece them together.

%--------------------------
\underline{(i)\ 
 Estimates for $\|u(t)\|_{L^2}$ and $\|v(t)\|_{L^2}$.
}
We first remark that \eqref{est_before} and  Lemma~\ref{lemma_asympt}
yield
\[
 \|u(t)\|_{W^{1,\infty}} + \|v(t)\|_{W^{1,\infty}}
 \le
 \frac{C\eps^{2/3}}{(1+t)^{1/2}}.
\]
Then we see from the standard energy method that  
\begin{align*}
 \frac{d}{dt}\Bigl(\|u(t)\|_{L^2} + \|v(t)\|_{L^2}\Bigr) 
 &\le 
 3|\kappa| \|v\|_{L^{\infty}}^2 \|\pa_x v\|_{L^2}
 +3|\lambda| \|u\|_{L^{\infty}}^2 \|\pa_x u\|_{L^2}
 \\
 &\le
 \frac{C\eps^2}{(1+t)^{1-\gamma/3}},
\end{align*}
whence
\begin{align*}
 \|u(t)\|_{L^2} + \|v(t)\|_{L^2}
 \le 
 C\eps 
 +
 \int_0^t \frac{C\eps^2}{(1+t')^{1-\gamma/3}} dt'
 \le
 C\eps(1+t)^{\gamma/3}.
\end{align*}

%--------------------------
\underline{(ii)\ 
 Estimates for $\|\op{J}_m u(t)\|_{L^2}$ and $\|\op{J}_{\mu} v(t)\|_{L^2}$.
} 
We note that \eqref{com_rel}, \eqref{id_J_and_P} and \eqref{nls_3}
yield
\[
 \op{J}_m\pa_x (v^3)=\op{P}(v^3) +2it\op{L}_{m}(v^3)
=3v^2\op{P}v +2v^3 +\op{L}_{m}(2itv^3)
\]
and
\[
 \op{P}v
=
\op{J}_{\mu}\pa_x v-2it\op{L}_{\mu}v
=
-v+\pa_x \op{J}_{\mu}v -2i\lambda t \pa_x (u^3).
\]
We also remember the commutation relation $[\op{L}_m,\op{J}_m]=0$. 
From them it follows that
\[
 \op{L}_m\bigl( \op{J}_mu-2i\kappa tv^3\bigr)=r,
\]
where
\[
r=
 -\kappa v^3 + 3\kappa v^2\pa_x \op{J}_{\mu}v
-6i \kappa \lambda t v^2 u^2\pa_x u.
\]
Since this $r$ can be estimated as 
\[
 \|r\|_{L^2} \le \frac{C\eps^2}{(1+t)^{1-\gamma/3}},
\]
the standard energy method leads to 
\[
\|\op{J}_mu-2i\kappa tv^3\|_{L^2}
\le
 \|xu_0\|_{L^2} + \int_0^t \frac{C\eps^2}{(1+t')^{1-\gamma/3}} dt'
 \le
 C\eps(1+t)^{\gamma/3}.
\]
Therefore we obtain 
\[
 \|\op{J}_m u(t)\|_{L^2}
 \le 
 2 t|\kappa|\|v\|_{L^{\infty}}^2\|v\|_{L^2} +C\eps(1+t)^{\gamma/3}
 \le 
 C\eps (1+t)^{\gamma/3}.
\]
In the same way, we have 
\[
 \|\op{J}_{\mu} v(t)\|_{L^2}
 \le 
 C\eps (1+t)^{\gamma/3}.
\]

%--------------------------
\underline{(iii)\ 
 Estimates for $\|\pa_x^2 u(t)\|_{L^2}+\|\pa_x^2 v(t)\|_{L^2}$.
} 
We  apply Lemma \ref{lemma_smoothing} with 
$\eta=\eps^{-2/3}$, $(\mu_1,\mu_2)=(m,\mu)$, $\mathbf{w}=(u,v)$ 
and $f=\pa_x^2 u$ or $\pa_x^2 v$. Then we obtain
\begin{align*}
 &\frac{d}{dt} \|S_{\Phi,m} \pa_x^2 u(t)\|_{L^2}^2
 + 
 \frac{1}{|m|}
 \int_{\R} \Phi(t,x) \Bigl| S_{\Phi,m} |\pa_x|^{1/2} \pa_x^2 u(t,x)\Bigr|^2dx
 \\
&\le
 2|\kappa| \left|
 \Jb{S_{\phi,m} \pa_x^2 u, S_{\Phi,m} \pa_x^{3} (v^3)}_{L^2}
 \right|
 + 
CB(t) \|\pa_x^2 u(t)\|_{L^2}^2
\end{align*}
as well as
\begin{align*}
 &\frac{d}{dt} \|S_{\Phi,\mu} \pa_x^2 v(t)\|_{L^2}^2
 + 
 \frac{1}{|m|}
 \int_{\R} \Phi(t,x) \Bigl| S_{\Phi,\mu} |\pa_x|^{1/2} \pa_x^2 v(t,x)\Bigr|^2dx
 \\
&\le
 2|\lambda| \left|
 \Jb{S_{\phi,\mu} \pa_x^2 v, S_{\Phi,\mu} \pa_x^{3} (u^3)}_{L^2}
 \right|
 + 
CB(t) \|\pa_x^2 v(t)\|_{L^2}^2,
\end{align*}
where 
\begin{align*}
 B(t)
 &= 
 e^{\frac{C}{\eps^{2/3}} ( \|u\|_{L^2}^2+\|v\|_{L^2}^2)}
 \Bigl(
  \eps^{-2/3}\|u\|_{W^{1,\infty}}^2  
 +
 \eps^{-2/3}\|v\|_{W^{1,\infty}}^2
   + \eps^{-2} \|u\|_{L^{\infty}}^6\\
 &\hspace{8mm}  + \eps^{-2} \|v\|_{L^{\infty}}^6
 +
  \eps^{-2/3}\|u\|_{L^2} \|\kappa \pa_x (v^3)\|_{L^2}
 +
  \eps^{-2/3}\|v\|_{L^2} \|\lambda \pa_x (u^3)\|_{L^2}
 \Bigr).
\end{align*}
Since \eqref{est_before} yields a rough $L^2$-bound
\[
 \|u(t)\|_{L^2}^2 + \|v(t)\|_{L^2}^2 
 \le
 C\Bigl(\int_{\R} \frac{d\xi}{\jb{\xi}^2}\Bigr)
 \sup_{\xi\in \R} 
 \Bigl(\jb{\xi} \bigl(|\alpha(t,\xi)|+|\beta(t,\xi)| \bigr)\Bigr)^2
 \le
 C\eps^{4/3},
\]
we see that $B(t)$ can be dominated by $C\eps^{\frac{2}{3}} (1+t)^{-1}$. 
We also observe that 
\[
 \pa_x^3 (v^3)=3v^2\pa_x(\pa_x^2v)+\rho_{1}
\]
with a remainder term $\rho_{1}$ satisfying 
\[
 \|\rho_{1}\|_{L^2}
 \le 
 C\|v\|_{W^{1,\infty}}^2\|v\|_{H^2}
 \le
 \frac{C\eps^{2}}{(1+t)^{1-\gamma/3}}.
\]
So it follows from Lemma~\ref{lemma_aux} that 
\begin{align*}
 &\left|\Jb{S_{\Phi,m}\pa_x^2 u, S_{\Phi,m} \pa_x^3(v^3)}_{L^2} \right|
 \\
 &\le
 3\left|\Jb{S_{\Phi,m}\pa_x^2 u, S_{\Phi,m} v^2 \pa_x(\pa_x^2 v)}_{L^2} \right|
 + \|S_{\Phi,m} \pa_x^2 u\|_{L^2} \|S_{\Phi,m} \rho_{1}\|_{L^2}\\
 &\le
 C_1\eps^{2/3} 
  \Bigl(
 \int_{\R} \Phi(t,x) \Bigl| S_{\Phi,m}|\pa_x|^{1/2} \pa_x^2 u(t,x) \Bigr|^2 dx
+
  \int_{\R} \Phi(t,x) \Bigl| S_{\Phi,\mu}|\pa_x|^{1/2} \pa_x^2 v(t,x) \Bigr|^2 
 dx
 \Bigr)\\
 &\hspace{6mm}+
  \frac{C\eps^{8/3}}{(1+t)^{1-2\gamma/3}}
\end{align*}
with some positive constant $C_1$ not depending on $\eps$. 
Similarly we have 
\begin{align*}
 &\left|\Jb{S_{\Phi,\mu}\pa_x^2 v, S_{\Phi,\mu} \pa_x^3(u^3)}_{L^2} \right|\\
 &\le
 C_1\eps^{2/3} 
  \Bigl(
 \int_{\R} \Phi(t,x) \Bigl| S_{\Phi,m}|\pa_x|^{1/2} \pa_x^2 u(t,x) \Bigr|^2 dx
+
  \int_{\R} \Phi(t,x) \Bigl| S_{\Phi,\mu}|\pa_x|^{1/2} \pa_x^2 v(t,x) \Bigr|^2 
 dx
 \Bigr)\\
 &\hspace{6mm}+
  \frac{C\eps^{8/3}}{(1+t)^{1-2\gamma/3}}.
\end{align*}
Summing up, we obtain 
\begin{align*}
 \frac{d}{dt}\Bigl(
 \|S_{\Phi,m}\pa_x^2 u(t)\|_{L^2}^2 + \|S_{\Phi,m}\pa_x^2 u(t)\|_{L^2}^2
 \Bigr)
 \le
  \frac{C\eps^{8/3}}{(1+t)^{1-2\gamma/3}}
\end{align*}
if $\eps$ is so small to satisfy 
\[
 2C_1(|\kappa|+|\lambda|)\eps^{2/3} 
 \le 
 \frac{1}{|m|} + \frac{1}{|\mu|}.
\]
By integration with respect to $t$, we obtain 
\[
  \|S_{\Phi,m}\pa_x^2 u(t)\|_{L^2}^2 + \|S_{\Phi,\mu}\pa_x^2 v(t)\|_{L^2}^2
 \le 
 C\eps^2 +C\eps^{2}(1+t)^{2\gamma/3},
\]
whence
\[
  \|\pa_x^2 u(t)\|_{L^2} + \|\pa_x^2 v(t)\|_{L^2}
 \le 
 C\eps(1+t)^{\gamma/3}.
\]

%--------------------------
\underline{(iv)\ 
 Estimates for $\|\op{P}u(t)\|_{L^2}+\|\op{P} v(t)\|_{L^2}$.
} 
By the commutation relations 
$[\op{L}_m, \op{P}]=2\op{L}_m$, $[\pa_x, \op{P}]=\pa_x$ 
and the Leibniz rule for $\op{P}$, we have 
\[
\op{L}_m(\op{P} u)
= (\op{P}+2)\bigl(3\kappa v^2\pa_x v\bigr)
=3\kappa v^2 \pa_x(\op{P} v)+\rho_2
\]
with a remainder term $\rho_2$ satisfying 
\[
 \|\rho_2\|_{L^2}
 \le 
 C \|v\|_{W^{1,\infty}}^2 
 \bigl( \|v\|_{H^1} + \|\op{P}v\|_{L^2}\bigr).
\]
Also we have
\[
\op{L}_{\mu}(\op{P} v)
=3\lambda u^2 \pa_x(\op{P} u)+\rho_3
\]
with 
\[
 \|\rho_3\|_{L^2}
 \le 
 C \|u\|_{W^{1,\infty}}^2 
 \bigl( \|u\|_{H^1} + \|\op{P}u\|_{L^2}\bigr).
\]
Thus, we can deduce as in the previous case that 
\[
  \|S_{\Phi,m}\op{P} u(t)\|_{L^2}^2 + \|S_{\Phi,\mu}\op{P} v(t)\|_{L^2}^2
 \le 
 C\eps^{2}(1+t)^{2\gamma/3},
\]
whence
\[
  \|\op{P} u(t)\|_{L^2} + \|\op{P} v(t)\|_{L^2}
 \le 
 C\eps(1+t)^{\gamma/3}.
\]

%--------------------------
\underline{(v)\  Conclusion.} 
Let $\phi=u$ or $v$. We note that the commutation relations \eqref{com_rel} give us 
\[
 \|\phi\|_{H^2}+\|\op{J}_m \phi\|_{H^1}
 \le 
C\bigl(
 \|\phi\|_{L^2}+ \|\pa_x^2 \phi\|_{L^2}
 +\|\op{J}_m \phi\|_{L^2}+ \|\op{J}_m \pa_x \phi\|_{L^2}\bigr)
\]
and that \eqref{id_J_and_P} yields
\[
 \|\op{J}_m \pa_x \phi\|_{L^2}
 \le
 \|\op{P} \phi\|_{L^2} + 2t\|\op{L}_m\phi\|_{L^2}.
\]
Therefore, by piecing together the estimates obtained in (i), (ii), (iii) 
and (iv), we arrive at the desired estimate \eqref{est_L2}.\qed

%----------------------------------------------------%
\subsection{ Estimates for $\alpha(t,\xi)$ and $\beta(t,\xi)$} 
\label{subsec_est_alpha}
%----------------------------------------------------%
The goal of this part is to prove \eqref{est_ptwise}. 
When $0\leq t\le 1$, the desired estimates follows immediately from the Sobolev 
embedding. Hence we have only to consider the case of $t\in [1,T)$. 

It follows from the definition of $\alpha$ that 
\begin{align}
\label{diff_alpha}
 i\pa_t \alpha(t,\xi)
 &=
 \frac{1}{1+im\xi} \op{F}_m \op{U}_{m}(t)^{-1}
 \Bigl[(1+\pa_x)\op{L}_m u(t,\cdot)\Bigr](\xi)\\
 &=
 \frac{\kappa}{1+im\xi}
 \Biggl(
 \op{F}_m \op{U}_{m}(t)^{-1}\Bigl[\pa_x(v^3)\Bigr] 
 +
 \op{F}_m \op{U}_{m}(t)^{-1}\Bigl[\pa_x^2(v^3)\Bigr]
 \Biggr).
 \nonumber
\end{align}
Next we put $\beta^{(l)}(t,\xi)=(i\mu \xi)^{l}\beta(t,\xi)$ 
so that 
$\pa_x^lv=\op{M}_{\mu}(t)\op{D}(t)\op{W}_{\mu}(t)\beta^{(l)}$ 
for $l=0, 1$. We also set
$\bigl(\op{E}^{\omega}(t)f\bigr)(y)=e^{i\frac{\omega}{2} ty^2} f(y)$ 
for $\omega \in \R$. Then we have 
\begin{align*}
\pa_x (v^3)
&=
3v^2\pa_x v\\
&=
3\Bigl(\op{M}_{\mu}(t) \op{D}(t) \op{W}_{\mu}(t) \beta \Bigr)^2 
\op{M}_{\mu}(t) \op{D}(t) \op{W}_{\mu}(t) \beta^{(1)}\\
&=
 \frac{3}{t}\op{M}_{3\mu}(t) \op{D}(t)
 \Bigl[
  \bigl( \op{W}_{\mu}(t) \beta \bigr)^2\op{W}_{\mu}(t) \beta^{(1)} 
 \Bigr].
\end{align*}
By the relation 
$\op{F}_m \op{U}_{m}(t)^{-1}\op{M}_{\omega}(t)\op{D}(t)
= \op{W}_{m}(t)^{-1}\op{E}^{\omega -m}(t)$, 
we see that 
\begin{align*}
 &\op{F}_m \op{U}_{m}(t)^{-1}\Bigl[\pa_x(v^3)\Bigr]\\
 %&=
 %\frac{3}{t}\op{W}_{m}(t)^{-1} \op{D}(t)^{-1}\op{M}_{3\mu-m}(t)
 %\op{D}(t)
 %\Bigl[
 % \bigl( \op{W}_{\mu}(t) \beta \bigr)^2\op{W}_{\mu}(t) \beta^{(1)} 
 %\Bigr]\\
 &=
  \frac{3}{t}\op{W}_{m}(t)^{-1} \op{E}^{3\mu-m}(t)
 \Bigl[
  \bigl( \op{W}_{\mu}(t) \beta \bigr)^2\op{W}_{\mu}(t) \beta^{(1)} 
 \Bigr]\\
 &=
  \frac{3}{t}\op{W}_{m}(t)^{-1} \op{E}^{3\mu-m}(t)
 \Bigl[i\mu\xi \beta^3 \Bigr]\\
 &\hspace{4mm}+
  \frac{3}{t}\op{W}_{m}(t)^{-1} \op{E}^{3\mu-m}(t)
 \Bigl[
  \bigl( \op{W}_{\mu}(t) \beta \bigr)^2\op{W}_{\mu}(t) \beta^{(1)} 
 -
  \beta^2 \beta^{(1)}
 \Bigr].
\end{align*}
Now we are interested in the principal part of the first term. 
Because of the relation 
$\op{W}_m(t)^{-1}=\op{D}(1/m) \op{U}_1(m/t) \op{D}(m)$, we have 
\begin{align*}
\Bigl(\op{W}_{m}(t)^{-1} \op{E}^{3\mu-m}(t) f\Bigr)(\xi)
 &=
 \sqrt{\frac{t}{2\pi mi}}\int_{\R}
 e^{i\frac{t}{2m}(m\xi-y)^2} e^{i\frac{t}{2}(3\mu-m) (\frac{y}{m})^2} 
  f\bigl(\frac{y}{m} \bigr)\, dy\\
 &=
 \sqrt{\frac{mt}{2\pi i}}\int_{\R}
 e^{i\frac{t}{2}(3\mu z^2 -2m\xi z +m\xi^2)} 
  f(z)\, dz\\
 &=
 e^{i\frac{m(3\mu-m)}{6\mu}t\xi^2}\sqrt{\frac{mt}{2\pi i}}\int_{\R}
  e^{i\frac{3\mu t}{2}(z-\frac{m\xi}{3\mu})^2} f(z)\, dz\\
 &=
 \gamma_1 \op{E}^{2\omega_1}(t) \op{D}(\frac{3\mu}{m}) \op{W}_1(3\mu t)^{-1}f (\xi)\\
 &=
 \gamma_1\sqrt{\frac{m}{3\mu}} e^{i\omega_1 t\xi^2}
 f \bigl(\frac{m\xi}{3\mu} \bigr)
 +
 \gamma_1 \op{E}^{2\omega_1}(t) 
 \op{D}(\frac{3\mu}{m}) \bigl(\op{W}_1(3\mu t)^{-1}-1\bigr)f (\xi)
\end{align*}
with an appropriate constant $\gamma_1 \in \C$  and 
\[
\omega_1=\frac{m(3\mu-m)}{6\mu}.
\]
So we can split 
$\op{F}_m \op{U}_{m}(t)^{-1}\Bigl[\pa_x(v^3)\Bigr]$ as 
\begin{align}
\label{decomp_N1}
 \op{F}_m \op{U}_{m}(t)^{-1}\Bigl[\pa_x(v^3)\Bigr]
 &=
\frac{\xi}{t} e^{i \omega_1 t\xi^2}  \Lambda_1(t,\xi) 
+\frac{\sigma_{1}(t,\xi)}{t},
\end{align}
where 
\[
 \Lambda_1(t,\xi)=i \gamma_1 m\sqrt{\frac{m}{3\mu}}  
            \biggl(\beta \bigl(t,\frac{m\xi}{3\mu}\bigr)\biggr)^3
\]
and
\begin{align*}
 \sigma_1(t,\xi)
 =&
 3\op{W}_{m}(t)^{-1} \op{E}^{3\mu-m}(t)
 \Bigl[
  \bigl( \op{W}_{\mu}(t) \beta \bigr)^2\op{W}_{\mu}(t) \beta^{(1)} 
 -
  \beta^2 \beta^{(1)}
 \Bigr]\\
 &+
 3\gamma_1 \op{E}^{2\omega_1}(t) 
 \op{D}(\frac{3\mu}{m}) \bigl(\op{W}_1(3\mu t)^{-1}-1\bigr)\Bigl[\beta^2 \beta^{(1)} 
 \Bigr].
\end{align*}
By virtue of the inequality
\[
 \bigl\| (\op{W}_m(t)-1)f \bigr\|_{L^{\infty}}
+
 \bigl\| (\op{W}_m(t)^{-1}-1)f \bigr\|_{L^{\infty}}
\le Ct^{-1/4} \|f\|_{H^1},
\]
we can see that 
\begin{align}
 |\sigma_{1}(t,\xi)| \le \frac{C\eps^2}{t^{1/4 -\gamma}}. 
\label{est_sigma}
\end{align}
Next we focus on the second term in \eqref{diff_alpha}. 
We first observe that 
\begin{align}
\label{decomp_N2}
 \pa_x^2(v^3)
 &=
 9v(\pa_x v)^2
 +\frac{3\mu}{it}
 \biggl( v^2 (\op{J}_{\mu}\pa_x v) -
         v(\pa_x v) (\op{J}_{\mu} v)\biggr)\\
 &=
 \frac{9}{t}\op{M}_{3\mu}(t) \op{D}(t)
 \Bigl[
   \op{W}_{\mu}(t) \beta \bigl( \op{W}_{\mu}(t) \beta^{(1)} \bigr)^2
 \Bigr]
 +
 \frac{3\mu}{it}
 \biggl( v^2 (\op{J}_{\mu}\pa_x v) -
         v(\pa_x v) (\op{J}_{\mu} v)\biggr).
\nonumber
\end{align}
By Lemma~\ref{lemma_prod}, we also have 
\begin{align*}
 \Bigl\|\op{F}_m \op{U}_{m}(t)^{-1}\Bigl[v^2 (\op{J}_{\mu}\pa_x v) -
         v(\pa_x v) (\op{J}_{\mu} v)\Bigr]\Bigr\|_{L^{\infty}}
 &\le
 \frac{C}{t^{1/4}}(C\eps^{2/3} t^{\gamma/3})^3
 =
 \frac{C\eps^2}{t^{1/4-\gamma}}.
\end{align*}
Therefore we can deduce that 
\begin{align*}
 \op{F}_m \op{U}_{m}(t)^{-1}\Bigl[\pa_x^2 (v^3)\Bigr]
 &=
 \frac{9}{t}\op{W}_{m}(t)^{-1} \op{E}^{3\mu-m}(t)
 \Bigl[(i\mu\xi)^2 \beta^3 \Bigr]
 +O(\eps^2t^{-5/4+\gamma})\\
 &=
 \frac{9 \gamma_1}{t}\sqrt{\frac{m}{3\mu}} e^{i\omega_1 t\xi^2} 
 \Bigl(i\mu\frac{m\xi}{3\mu} \Bigr)^2 
 \Bigl(\beta \bigl(t,\frac{m\xi}{3\mu} \bigr) \Bigr)^3
+O(\eps^2t^{-5/4+\gamma})\\
 &=
i m\xi \frac{\xi}{t} e^{i\omega_1 t\xi^2} \Lambda_1(t,\xi)  
+O(\eps^2t^{-5/4+\gamma}).
\end{align*}
By \eqref{diff_alpha}, \eqref{decomp_N1}, \eqref{est_sigma} and 
\eqref{decomp_N2}, we have
\begin{align}
 i\pa_t \alpha
 &=
\frac{\xi}{t} e^{i \omega_1 t\xi^2} A_1(t,\xi) +R_{1}(t,\xi)
\label{eq_reduced}
\end{align}
with 
\[
 A_1(t,\xi)=\frac{\kappa (1+im\xi)}{1+im\xi} \Lambda_1(t,\xi)
=\kappa \Lambda_1(t,\xi)\]
and a remainder term $R_1(t,\xi)$ satisfying 
\[
 |R_{1}(t,\xi)| \le \frac{C\eps^2}{\jb{\xi} t^{5/4 -\gamma}}. 
\]
By using the identity
\[
 \frac{\xi  e^{i\omega t\xi^2}}{t} A(t,\xi)
=
 i \pa_t 
 \left(\frac{ -i\xi e^{i\omega t\xi^2}}{1+i \omega t\xi^2} A(t,\xi)\right) 
 -
 t e^{i \omega t\xi^2}
 \pa_t\left( \frac{ \xi A(t,\xi)}{t(1+i\omega t\xi^2)}  \right)
\]
as well as the inequality
\[
 \sup_{\xi \in \R} 
 \left|\frac{\xi }{1+i \omega t\xi^2} \right| 
 \le 
 \frac{C}{t^{1/2}},
\]
we see that the first term in the right-hand side of \eqref{eq_reduced} 
can be splitted into the following form:
\[
 \frac{\xi e^{i\omega_{1}t\xi^{2}}}{t}A_1 =i\pa_t A_{1,1} +A_{1,2}, \quad 
 |A_{1,1}(t,\xi)| \le \frac{C\eps^2}{t^{1/2} \jb{\xi}^3}, \quad
 |A_{1,2}(t,\xi)| \le \frac{C\eps^2}{t^{3/2} \jb{\xi}}.
\]
Therefore we have 
\begin{align*}
 |\alpha(t,\xi)-A_{1,1}(t,\xi)|
 &\le 
 |\alpha(t,1)-A_{1,1}(t,1)|+ \int_1^{t} |A_{1,2}(t',\xi)+R_1(t',\xi)|dt'
 \le 
 \frac{C\eps}{\jb{\xi}},
\end{align*}
whence
\[
 \jb{\xi}|\alpha(t,\xi)|
 \le 
 \jb{\xi}|A_{1,1}(t,\xi)|+ \jb{\xi}|\alpha(t,\xi)-A_{1,1}(t,\xi)|
 %+
 %\jb{\xi} \int_1^{t} |A_{1,2}(t',\xi)+R_1(t',\xi)|dt'
 \le C\eps.
\]
In the same way, we have 
\begin{align*}
 i\pa_t \beta(t,\xi)
 &=
 \frac{\xi}{t} e^{i\omega_2 t\xi^2} A_2(t,\xi) +R_{2}(t,\xi),
\end{align*}
where 
\[
 \omega_2=\frac{\mu (3m-\mu)}{6m},
\]
\[
 A_2(t,\xi)=i\lambda \gamma_2 \mu \sqrt{\frac{\mu }{3m}} 
  \biggl(\alpha\bigl(t,\frac{\mu \xi}{3m}\bigr)\biggr)^3,
\]
\[
 |R_{2}(t,\xi)| \le \frac{C\eps^2}{\jb{\xi} t^{5/4 -\gamma}},
\] 
and  $\gamma_2 \in \C$ is an appropriate constant. 
Therefore we can deduce as before that 
\[
 \jb{\xi}|\beta(t,\xi)|\le C\eps.
\]
Summing up, we obtain the desired estimate \eqref{est_ptwise}.
\qed

%---------------------------------------------------------------------------%
\section{Proof of Theorem~\ref{thm_main}} \label{Sec_proof1}
%---------------------------------------------------------------------------%
In this section, we will prove Theorem~\ref{thm_main}. 
We first recall the local existence. 
For fixed  $t_0\ge 0$, let us consider the initial value 
problem 
\begin{align}
\left\{\begin{array}{lc}
 \begin{array}{l}
 \op{L}_{m} u = \kappa \pa_x (v^3), \\
 \op{L}_{\mu} v =\lambda \pa_x (u^3), 
 \end{array}
 & t>t_0,\ x\in \R,\\
(u,v)\bigm|_{t=t_0}=(\varphi,\psi),
& x \in \R.
\end{array}\right.
\label{ivp_shift}
\end{align}

%----------------------
\begin{lem} \label{lemma_local}
Let  $\varphi$, $\psi \in H^2\cap H^{1,1}$. 
There exists a positive constant $\eps_0$, not depending on $t_0$, 
such that the following holds: 
for any $\underline{\eps} \in (0,\eps_0)$ and $M\in (0,\infty)$, 
one can choose a positive constant $\tau^*=\tau^{*}(\underline{\eps},M)$, 
which is independent of $t_0$, 
such that \eqref{ivp_shift}--\eqref{data} admits a unique pair of solutions 
$u$, $v \in C([t_0,t_0+\tau^*]; H^2\cap H^{1,1})$, 
provided that $\|\varphi\|+\|\psi\|_{L^2}  \le \underline{\eps}$ 
and
\[
 \sum_{l=0}^{1} \Bigl(
 \Bigl\|\bigl(x+\frac{it_0}{m}\pa_x \bigr)^l \varphi \Bigr\|_{H^{2-l}} 
+
 \Bigl\|\bigl(x+\frac{it_0}{\mu}\pa_x \bigr)^l\psi \Bigr\|_{H^{2-l}} 
 \Bigr)\le M. 
\]
\end{lem}
%-----------------------
We omit the proof of this lemma because it is standard.
\\

Now we are going to prove the global existence by the so-called 
bootstrap argument. 
Let $T^*$ be the supremum of all $T \in (0,\infty]$ such that 
the problem \eqref{nls_3}--\eqref{data} admits a unique pair of solutions 
$u$, $v \in C([0,T);H^2\cap H^{1,1})$. 
By Lemma \ref{lemma_local} with $t_0=0$, we have $T^*>0$ if 
$\|u_0\|_{L^2}+\|v_0\|_{L^2}\le \eps< \eps_0$. 
We also set 
\[
 T_{*}
=\sup \bigl\{ \tau \in [0,T^*)\, |\, E(\tau) \le \eps^{2/3} \bigr\}.
\]
Note that $T_*>0$ because of the continuity of 
$[0,T^*) \ni \tau \mapsto E(\tau)$ and 
$E(0) \le C\eps \le \frac{1}{2}\eps^{2/3}$ 
if $\eps$ is suitably small. 
We claim that $T_*=T^{*}$ if $\eps$ is small enough. Indeed, if 
$T_{*}<T^*$, Lemma \ref{lemma_apriori} with $T=T_*$ yields 
\[
 E(T_*) \le K\eps \le \frac{1}{2}\eps^{2/3}
\] 
for $\eps\le \eps_2:=\min\{\eps_1, 1/(2K)^3\}$, 
where $K$ and $\eps_1$ are mentioned in Lemma \ref{lemma_apriori}. 
By the continuity of $[0,T^*)\ni \tau \mapsto E(\tau) $, 
we can take $T^{\flat} \in (T_*,T^*)$ such that 
$E(T^{\flat}) \le \eps^{2/3}$, which contradicts the 
definition of $T_*$. Therefore we must have $T_*=T^*$. 
By using Lemma \ref{lemma_apriori} with $T=T^*$ again, we see that 
\[
  \sum_{l=0}^{1}\Bigl(\|J_{m}^{l} u(t,\cdot)\|_{H^{2-l}}
  +
  \|J_{\mu}^{l} v(t,\cdot)\|_{H^{2-l}} \Bigr)
  \le 
  K\eps (1+t)^{\frac{\gamma}{3}}, 
\]
and
\[
 \sup_{\xi \in \R} \jb{\xi} \Bigl(|\alpha(t,\xi)|+|\beta(t,\xi)|\Bigr)
 \le K\eps
\]
for $t\in [0,T^*)$ and $\eps\leq \eps_{1}$. In particular we have
\[
 \sup_{t\in [0,T^*)}\Bigl(\|u(t)\|_{L^2} +\|v(t)\|_{L^2}\Bigr)
 \le 
  C\sup_{(t,\xi)\in [0,T^*) \times \R} 
 \jb{\xi} \Bigl(|\alpha(t,\xi)| + |\beta(t,\xi)| \Bigr)
 \le 
 C^{\flat}\eps
\]
for $\eps\leq \eps_{1}$ with some $C^{\flat}>0$. 
Next we assume $T^*<\infty$. Then, by setting 
$\eps_3=\min\{\eps_2, \eps_0/2C^{\flat}\}$ 
and $M=K\eps_3 (1+T^*)^{\gamma/3}$, we have 
\[
  \sup_{t\in [0,T^*)} \sum_{l=0}^{1} \Bigl( 
 \|J_{m}^{l} u(t,\cdot)\|_{H^{2-l}} 
 +
 \|J_{\mu}^{l} v(t,\cdot)\|_{H^{2-l}} 
 \Bigr)\le M
\]
as well as 
\[
 \sup_{t\in [0,T^*)}
  \bigl(\|u(t)\|_{L^2}+\|v(t)\|_{L^2}\bigr) \le \eps_0/2<\eps_0
\]
for $\eps \le \eps_3$. By Lemma \ref{lemma_local}, 
there exists $\tau^*>0$ such that \eqref{nls_3}--\eqref{data} admits a 
unique pair of solutions $u, v \in C([0, T^*+\tau^*);H^{2}\cap H^{1,1})$. 
This contradicts the definition of $T^*$, which means $T^*=\infty$ for 
$\eps \in (0,\eps_3]$. 
This completes the proof of the global existence part of Theorem \ref{thm_main}. We also conclude that the estimates \eqref{est_L2} and \eqref{est_ptwise} are 
valid for $(t,\xi) \in [0,\infty)\times \R$.

Next we turn our attentions to the asymptotic behavior.
For given $\delta>0$, we set $\gamma=\min\{\delta, 1/5\} \in (0,1/4)$. 
Recalling the argument in \S \ref{subsec_est_alpha}, we see that
\[
 i\pa_t\Bigl(\alpha(t,\xi) -A_{1,1}(t,\xi)\Bigr)
 = A_{1,2}(t,\xi)+R_1(t,\xi)
\]
with 
\[
 |A_{1,1}(t,\xi)|\le \frac{C\eps^2}{\jb{\xi}^3 t^{1/2}},
 \quad
 |A_{1,2}(t,\xi)|\le \frac{C\eps^2}{\jb{\xi}t^{3/2}}, 
 \quad
 |R_1(t,\xi)| 
 \le 
 \frac{C \eps^2}{\jb{\xi} t^{5/4-\gamma}}
\]
for $t\ge 1$, $\xi \in \R$. 
These estimates allow us to define 
$\alpha_+ \in L^2\cap L^{\infty}$ by 
\[
 \alpha_+(\xi):=\alpha(1,\xi) -A_{1,1}(1,\xi)
-i \int_1^{\infty} A_{1,2}(t',\xi) + R_1(t',\xi) dt'.
\]
Also we set 
$u_+:=\mathcal{F}_{m}^{-1} \alpha_+$. 
Then, because of the relation
\begin{align*}
 \alpha(t,\xi)-A_{1,1}(t,\xi)
 &= 
 \alpha(1,\xi) -A_{1,1}(1,\xi) 
 -i \int_1^{t}A_{1,2}(t',\xi)+R_1(t',\xi) dt'\\
 &=
  \alpha_+(\xi)+i \int_{t}^{\infty}A_{1,2}(t',\xi)+R_1(t',\xi) dt',
\end{align*}
we have
\begin{align*}
 \|\alpha(t) -\alpha_+ \|_{L^{2}\cap L^{\infty}}
 &\le 
 \|A_{1,1}(t)\|_{L^{2}\cap L^{\infty}}
 +
 \int_t^{\infty}  
  \|A_{1,2}(t')+R_1(t') \|_{L^{2}\cap L^{\infty}} 
  dt'
 \le 
 C\eps^2 t^{-1/4+{\gamma}}.
\end{align*}
Combining the result obtained above and using 
Lemma \ref{lemma_asympt}, we have
\begin{align*}
 \|u(t) -\op{U}_{m}(t)u_+\|_{L^2}
 &=
 \|
 \op{F}_{m}\op{U}_{m}(t)^{-1}u(t) -\op{F}_{m}u_+
 \|_{L^2}\\
 &=
 \|\alpha(t) -\alpha_+ \|_{L^2}\\
 &\le 
 C\eps^2 t^{-1/4+\gamma}
 \le 
 C\eps^2 t^{-1/4+\delta}
\end{align*}
and
\begin{align*}
 &\|
   u(t) -\op{M}_{m}(t) \op{D}(t)\op{F}_{m}u_+
 \|_{L^{\infty}}\\
 &\le
 \|
   u(t) -\op{M}_{m}(t) \op{D}(t)\op{F}_{m}\op{U}_{m}(t)^{-1}
  u(t)
 \|_{L^{\infty}} 
 +
  \|
   \op{M}_{m}(t) \op{D}(t)\bigl[\alpha(t)-\alpha_+\bigr]
 \|_{L^{\infty}}\\
 &\le
 Ct^{-3/4} (\|u(t)\|_{L^2}+\|\op{J}_{m}u(t)\|_{L^2} )
 + 
 Ct^{-1/2} \|\alpha(t)-\alpha_+\|_{L^{\infty}}\\
 &\le
 C\eps t^{-3/4+\gamma/3} + C\eps^2 t^{-1/2 -1/4+\gamma}\\
 &\le 
 C\eps t^{-3/4+\delta}
\end{align*}
for $t\geq 1$. Similarly, we can specify the large-time behavior 
of $v(t,x)$ with the aid of the asymptotics for $\beta(t,\xi)$.
\qed

%---------------------------------------------------------------------------%
\section{Concerning Theorem~\ref{thm_second}} \label{Sec_proof2}
In this section, we give an outline of the proof of Theorem~\ref{thm_second}. 
The argument is almost parallel to that for Theorem~\ref{thm_main}. 
So we only point out two main differences. 
%The other parts of the proof are exactly same. \\
\begin{description}
\item[(1)]\  Since the nonlinear term of the first equation in \eqref{nls_2} 
is not the divergence form, we need a modification to get the estimate for 
$\|\op{J}_m u\|_{L^2}$. If $\mu\ne 2m$, it holds that 
\[
 \cc{u}^2\pa_x v
=\frac{\mu}{\mu-2m}\pa_x(\cc{u}^2 v)
 +\frac{2m\mu}{it(\mu-2m)}\biggl(
 \cc{u}(\cc{\op{J}_m u})v - \cc{u}^2 \op{J}_{\mu} v 
\biggr).
\]
Therefore we can show that 
\[
 \biggl\|\op{J}_m u- \frac{2i\kappa \mu t}{\mu-2m}\cc{u}^2 v \biggl\|_{L^2}
 \le C\eps (1+t)^{\gamma/3}. 
\]
\item[(2)]\ 
The reduced equation satisfied by $\alpha$ becomes 
\[
i\pa_t\alpha(t,\xi)=\frac{\xi}{t}e^{i\omega_3 t\xi^2}A_3(t,\xi) +R_3(t,\xi)
\]
with
\[
 \omega_3=\frac{m(\mu-3m)}{2(\mu-2m)}.
\]
So, if $\mu\ne 3m$ and $\mu\ne 2m$, we can obtain the desired pointwise bound 
for $\alpha(t,\xi)$ by taking the oscillating factor into account.
\end{description}

%---------------------------
\begin{rmk} 
In the case of $\mu=2m$, both of the good things mentioned above are missing, 
and we have no idea how to treat this case.
\end{rmk}
%-----------------------------

%---------------------------
\begin{rmk} \label{rem_resonance}
The case $\mu=3m$ is covered by the previous work \cite{LS}, 
if we put the additional restrictions 
\begin{align}
\realpart(\kappa \lambda)<0 
\quad \mbox{and}\quad 
\imagpart(\kappa \lambda)=0
\label{assump_additional}
\end{align}
(cf.~the condition (b$_0$) in \cite{LS}). 
To see where these restrictions come from, let us focus on the reduced system 
satisfied by $(\alpha,\beta)$. 
In the case of $\mu=3m$, the argument analogous to \S~\ref{subsec_est_alpha} 
yields
\[
 \pa_t \alpha = \frac{\mu\xi \kappa}{t} \cc{\alpha}^2\beta+r_1,\qquad
 \pa_t \beta = \frac{3m \xi \lambda}{t} \alpha^3+r_2,
\]
where $r_1$, $r_2$ are suitable harmless terms. 
Under \eqref{assump_additional}, 
we can choose $c>0$ such that $\lambda=-c\cc{\kappa}$. With this $c$, we have 
\begin{align*}
 \pa_t\Bigl(c|\alpha|^2+|\beta|^2\Bigr)
 &=
 2\realpart\Bigl(
  c\cc{\alpha}\pa_t \alpha + \cc{\beta}\pa_t \beta
 \Bigr)\\
 &=
 \frac{2c\mu\xi}{t}\realpart\Bigl(
  \kappa \cc{\alpha}^3 \beta - \cc{\kappa} \alpha^3 \cc{\beta}
 \Bigr)+ \mbox{(harmless terms)}\\
 &= 0+\mbox{(harmless terms)}.
\end{align*}
This allows us to obtain the a priori estimate, and thus the small data global 
existence. However, the asymptotic profile of the solution 
cannot be specified in this case. 
This point should be contrasted with the non-resonance case. 
\end{rmk}
%-----------------------------

%---------------------------------------------------------------------------%
\section{A generalization} \label{Sec_gene}
In this section, we give a generalization of Theorems \ref{thm_main} and 
\ref{thm_second} to the $N$-component system \eqref{nls_N}. 
We set $I_N=\{1,\ldots, N\}$ and 
$\sh{I}_N =\{1,\ldots, N, N+1,\ldots, 2N \}$. 
For $\mathbf{z} =(z_j)_{j\in I_N}\in \C^N$, we write 
\[
 \sh{\mathbf{z}} =(\sh{z}_k)_{k \in \sh{I}_N}
 :=(z_1,\ldots,z_N, \overline{z_1},\ldots, \overline{z_N})
\in \C^{2N}.
\]
Then general cubic nonlinear term $\mathbf{F}=(F_j)_{j \in I_N}$ 
can be written as 
\begin{align*}
 F_j(\bfu,\pa_x \bfu)
 =\sum_{l_1, l_2,l_3=0}^{1}\sum_{k_1, k_2, k_3 \in \sh{I}_N}
 C_{j, k_1, k_2, k_3}^{l_1,l_2,l_3} (\pa_x^{l_1} \sh{u}_{k_1}) 
 (\pa_x^{l_2} \sh{u}_{k_2}) (\pa_x^{l_3} \sh{u}_{k_3})
\end{align*}
with suitable $C_{j, k_1, k_2, k_3}^{l_1,l_2,l_3} \in \C$. 
We also introduce the following notation: 
\[
 \tm_k=\left\{\begin{array}{cl}
 m_k& (k=1,\ldots, N),\\[4mm]
 -m_{(k-N)} & (k = N+1,\ldots, 2N).
 \end{array}\right.
\]
The following theorem is a natural generalization of 
Theorems \ref{thm_main} and \ref{thm_second}.

%-------------------------
\begin{thm} \label{thm_N}
Assume that the following {\rm (i)} and {\rm (ii)} are satisfied 
for all $j \in I_N$ and $k_1,k_2, k_3 \in \sh{I}_N$:  
%------------%
\begin{itemize}
\item[(i)] 
 $C_{j,k_1,k_2,k_3}^{0,0,0}=0$. 
\item[(ii)] $\tm_{k_1} + \tm_{k_2} + \tm_{k_3} \in \{0, m_j\}$ 
implies $C_{j,k_1,k_2,k_3}^{l_1,l_2,l_3}=0$ for $l_1, l_2, l_3 =0,1$.
\end{itemize}
%------------%
Let 
$\varphi=(\varphi_j)_{j \in I_ N} \in H^3\cap H^{2,1}$, 
and suppose that $\|\varphi\|_{H^3}+\|\varphi\|_{H^{2,1}}$ is sufficiently 
small. Then \eqref{nls_N} admits a unique global solution 
$\bfu=(u_j)_{j \in I_N} \in C([0,\infty); H^3\cap H^{2,1})$. 
Moreover, for each $j\in I_N$, there exists $\varphi_j^+ \in L^2(\R_{x})$ with 
$\hat{\varphi}_j^+ \in L^{\infty}(\R_{\xi})$ such that 
\[
 u_j(t)
 =
 e^{i\frac{t}{2m_j}\pa_x^2} \varphi_j^+ + O(t^{-1/4+\delta})
 \quad \mbox{in}\ L^2(\R_x)
\]
and 
\[
 u_j(t,x)=\sqrt{\frac{m_j}{i t}}\, \hat{\varphi}^+_j
 \left( \frac{m_j x}{t}\right) e^{i\frac{m_j x^2}{2t}} 
 + O(t^{-3/4+ \delta})
 \quad \mbox{in}\ L^{\infty}(\R_x) 
\]
as $t \to \infty$, where $\delta>0$ can be taken arbitrarily small.
\end{thm}
%-------------------------

%---------------------------
\begin{rmk} 
The assumption {\rm (i)} says that each $F_j$ contains at least one derivative 
of the unknowns. The assumption {\rm (ii)} means that the cubic interaction of 
the form 
$(\pa_x^{l_1} \sh{u}_{k_1})(\pa_x^{l_2} \sh{u}_{k_2})(\pa_x^{l_3} \sh{u}_{k_3})$ is permitted in the $j$-th equation only if 
$\tm_{k_1} + \tm_{k_2} + \tm_{k_3} \not\in \{0, m_j\}$. 
Note that the assumptions $\mu\ne 3m$ and $\mu\ne 2m$ in 
Theorem~\ref{thm_second} can be read as 
\[
(-m)+(-m)+\mu \not \in \{0,m\} \quad \mbox{and} \quad 
m+m+ m \not \in \{0,\mu\}. 
\]
\end{rmk}
%-----------------------------

Finally, let us give a sketch of the proof of Theorem~\ref{thm_N}. 
The main step of the proof is to get the a priori estimate for 
\[
 E(T) =\sup_{0\le t< T}
 \sum_{j\in I_N}\biggl[
 (1+t)^{-\frac{\gamma}{3}} \Bigl(
  \|u_j(t)\|_{H^{3}} + \|J_{m_j} u_j(t)\|_{H^{2}}
 \Bigr) 
 + 
 \sup_{\xi \in \R} \Bigl(\jb{\xi}^2 |\alpha_j(t,\xi)|\Bigr)
 \biggr],
\]
where
\[
\alpha_j(t,\xi)
 =
\op{F}_{m_j}\Bigl[\op{U}_{m_j}(t)^{-1}u_j(t,\cdot)\Bigr](\xi).
\] 
To control the $L^2$-norm of $\op{J}_{m_j} u_j$, we use the following 
algebraic lemma: 
%-------------------------
\begin{lem} \label{lem_alg}
If $\mu_1+\mu_2+\mu_3\ne 0$, we have
\[
 f_1f_2\pa_x f_3 =\frac{\mu_3}{\mu_1+\mu_2+\mu_3}\pa_x\Bigl(f_1f_2f_3\Bigr)
+\frac{R}{it(\mu_1+\mu_2+\mu_3)},
\]
where
\[
 R=\mu_2\mu_3f_1\Bigl( f_2\op{J}_{\mu_3}f_3 -(\op{J}_{\mu_2}f_2)f_3\Bigr)
+
\mu_1\mu_3f_2\Bigl(f_1\op{J}_{\mu_3}f_3 -(\op{J}_{\mu_1}f_1)f_3\Bigr).
\]
\end{lem}
%-------------------------
By virtue of this lemma, we can choose a suitable cubic term 
$\Gamma_j$ such that $\|\op{J}_{m_j} u_j(t) -\Gamma_j(t)\|_{L^2}$ 
is dominated by $C\eps(1+t)^{\gamma/3}$ 
under the assumption $E(T)\le \eps^{2/3}$. 
To obtain the pointwise estimate for $\alpha_j(t,\xi)$, 
we derive the reduced equation including the oscillating factor 
$\frac{\xi}{t}e^{i \omega t\xi^2}$ with
\[
\omega
=
\frac{m_j^2}{2}
 \biggl(
  \frac{1}{m_j} - \frac{1}{\tm_{k_1}+\tm_{k_2}+\tm_{k_3}} 
 \biggr)
\in \R\backslash\{0\}.
\]
This enables us to get an a priori bound for $\alpha_j(t,\xi)$. 
The other parts of the proof are essentially the same as those for 
Theorems~\ref{thm_main} and \ref{thm_second}, so we omit the detail.

%----------------------------------------------------------------------------%
\subsection*{Acknowledgments}

The work of C.~L. is supported by NNSFC under Grant No. 11461074. 
The work of H.~S. is supported by Grant-in-Aid for Scientific Research~(C) 
(No.~25400161), JSPS.

%%%%%%%%%%%%%%%%%%%%%%%%%%%%%%%%%
% References
%%%%%%%%%%%%%%%%%%%%%%%%%%%%%%%%%

%%%%%%%%%%%%%%%%%%%%%%%%%%%%%%%%%%%%%%%%%%%%%%%%%%%%%%%%%%%%%%%%%%%%%%%%%%%%%%%
\end{document}